\documentclass[11pt]{article}
\usepackage{amsmath, amsthm, amsfonts}
\usepackage[dvips]{graphics}
\usepackage[dvips]{graphicx}

\begin{document}

\pagestyle{plain} \setcounter{page}{1}

\noindent \textbf{{\LARGE Projected iterations of fixed point type to solve nonlinear partial Volterra integro--differential equations}}

\noindent M.I. Berenguer, D. G\'{a}mez

\noindent {\footnotesize \textsl{Universidad de Granada,E.T.S. Ingenier\'ia Edificaci\'{o}n, Departamento de
Matem\'atica Aplicada, c/ Severo Ochoa s/n, 18071 Granada (Spain). E-mails:
maribel@ugr.es, domingo@ugr.es}}


\begin{abstract}
{In this paper, we propose a method to approximate the fixed point of an operator in a Banach space.

Using biorthogonal systems, this method is applied to build an approximation of  the
solution of a class of nonlinear partial integro--differential equations. The theoretical findings are illustrated with several numerical examples, confirming the reliability, validity and precision of the proposed method. }


 
\end{abstract}

%
%
\linespread{2}

\section{Introduction}
\linespread{2}

The Fixed--Point Theory of operators is a major research area in nonlinear analysis and one of the most powerful and fruitful tools of modern mathematics. This flourishing area of research in pure and applied mathematics can be used to establish the existence and uniqueness of solutions to problems that arise naturally in applications. However the explicit calculation of the fixed point of an operator is only possible in some simple cases, and in most cases it is essential to approximate this fixed point by a computational method.

The study of nonlinear partial integro--differential equations (NPIDEs) has become a subject of considerable interest. These equations  and similar ones arise in a variety of science and technology fields such as reaction-diffusion problems (\cite{dixon}), theory of elasticity (\cite{silling}), heat conduction (\cite{alkhaled}), mechanic of solids (\cite{appell}),  population dynamics (\cite{shivanian}), transient, conductive and radiative transport (\cite{frankel}), and other applications. 

Various numerical methods have been developed  for approximating the solution of NPIDEs. For example in \cite{avazzadeh}, the method is based on radial basis functions and in \cite{bobo} the authors generalize the Lomov's regularization method. In \cite{gur} a Laguerre collocation method is presented to solve certain nonlinear partial integro--diferential equations,  in \cite{rezvan} on the fractional differential transform method, in \cite{singh} the numerical solution is computed using the 2D shifted orthogonal Legendre polynomial system and,  in \cite{hameed}  the Newton--Kantorovich method is used for similar equations.

Our purpose in this paper is: (i) to present a method to approximate the fixed point of a  operator defined in a Banach space, by means of a composition of operators defined in said space (ii)  develop and apply the method presented for get an approximation of the solution of an nonlinear partial Volterra integro--differential equation (NPVIDE).

Section 2 deals with point (i) and in Section 4 is developed (ii) and we include several examples that show the validity of the exposed method. In addition, in section 3 we briefly present certain basic notions of the theory of biorthogonal systems that are necessary in section 4.

 \section{Approximation of fixed points of operators in Banach spaces}

 The Banach Fixed--Point Theorem is well known
(\cite{jameson}): \textsl{Let
$(X,\|\cdot\|)$ be a Banach space, let $\textbf{F}:X \longrightarrow X$ and
let $\{\mu_n\}_{n \geq 1}$ be a sequence of nonnegative real
numbers such that the series $\sum_{n \ge 1}\mu_n$ is convergent
and for all $u_1, u_2 \in X$ and for all $n\geq 1$, $\|\textbf{F}^n u_1-\textbf{F}^n u_2\|\leq
\mu_n\|u_1-u_2\|$. Then $\textbf{F}$ has a unique fixed point $u\in X$. Moreover,
if $\widetilde{u}$ is an element in $X$, then we have that for all
$n\in \mathbb{N}$,
\begin{equation}\label{1}
\|\textbf{F}^n\widetilde{u}-u\|\leq
\left(\sum_{i=n}^{\infty}\mu_i \right) \|\textbf{F}\widetilde{u}-\widetilde{u}\|.
\end{equation}
In particular, }
\begin{equation} \label{limite}
 u=\lim_n \textbf{F}^n(\widetilde{u}) \end{equation} 
 
Note that (\ref{limite}) implies that, for $n$ large enough ,$\|\textbf{F}^n\widetilde{u}-u\|$  will be small enough.

 The difficulty in obtaining the fixed point of an operator as the previous limit is clear in practical situations where it is sometimes impossible to apply the operator $ \textbf{F} $.

The result that is expounded below ensures, under certain conditions, a way to approximate the fixed point $ u $ of the operator in question through a composition of other operators.

\noindent\textbf{Lemma 1.} \textsl{Let
$(X,\|\cdot\|)$ be a Banach space, let $\textbf{F}:X \longrightarrow X$ and
let $\mu_0=1$ and  $\{\mu_n\}_{n \geq 1}$ be a sequence of nonnegative real
numbers such that  for all $n\geq 1$, $\|\textbf{F}^n u_1-\textbf{F}^n u_2\|\leq
\mu_n\|u_1-u_2\|$. Let $\widetilde{u} \in X$ and  $\textbf{G}_1,\cdots, \textbf{G}_m:X \longrightarrow X $. Then 
 $$ \|\textbf{F}^m(\widetilde{u})-\textbf{G}_m\circ...\circ \textbf{G}_1(\widetilde{u})\|\leq \sum_{p=1}^{m} \mu_{m-p}  \|\textbf{F}\circ \textbf{G}_{p-1}\circ\cdots \circ \textbf{G}_1(\widetilde{u})-\textbf{G}_p\circ \cdots \circ \textbf{G}_1(\widetilde{u})\|.
 $$}
  
  \bigskip
  
\textsl{Proof.}   
$$ \|\textbf{F}^m(\widetilde{u})-\textbf{G}_m\circ...\circ \textbf{G}_1(\widetilde{u})\|\leq  $$
$$
\|\textbf{F}^m(\widetilde{u})-\textbf{F}^{m-1}\circ\textbf{G}_1 (\widetilde{u})\|+\|\textbf{F}^{m-1}\circ \textbf{G}_1 (\widetilde{u})-\textbf{F}^{m-2}\circ \textbf{G}_2 \circ \textbf{G}_1 (\widetilde{u})\|
+\dots+ $$ 
$$ \|\textbf{F}^{2}\circ \textbf{G}_{m-2}\circ\cdots \circ\textbf{G}_1(\widetilde{u})-\textbf{F}\circ\textbf{G}_{m-1}\circ \cdots \circ\textbf{G}_1(\widetilde{u})\|+ \|\textbf{F}\circ \textbf{G}_{m-1}\circ\cdots \circ \textbf{G}_1(\widetilde{u})-\textbf{G}_m \circ\cdots \circ \textbf{G}_1(\widetilde{u})\|= $$ 
$$
\left(\sum_{p=1}^{m-1}  \|\textbf{F}^{m-p+1}\circ \textbf{G}_{p-1}\circ\cdots \circ\textbf{G}_1(\widetilde{u})-\textbf{F}^{m-p}\circ\textbf{G}_p\circ \cdots \circ\textbf{G}_1(\widetilde{u})\|\right)+\|\textbf{F}\circ \textbf{G}_{m-1}\circ\cdots \circ \textbf{G}_1(\widetilde{u})-\textbf{G}_m\circ\cdots \circ \textbf{G}_1(\widetilde{u})\|= $$
$$\left(\sum_{p=1}^{m-1}  \|\textbf{F}^{m-p}\circ \textbf{F}\circ \textbf{G}_{p-1}\circ\cdots \circ\textbf{G}_1(\widetilde{u})-\textbf{F}^{m-p}\circ\textbf{G}_p\circ \cdots \circ\textbf{G}_1(\widetilde{u})\|\right) + \|\textbf{F}\circ \textbf{G}_{m-1}\circ\cdots \circ \textbf{G}_1(\widetilde{u})-\textbf{G}_m \circ\cdots \circ \textbf{G}_1(\widetilde{u})\|\leq $$
$$\left(\sum_{p=1}^{m-1} \mu_{m-p}  \|\textbf{F}\circ \textbf{G}_{p-1}\circ\cdots \circ \textbf{G}_1(\widetilde{u})-\textbf{G}_p\circ \textbf{G}_{p-1}\circ\cdots \circ \textbf{G}_1(\widetilde{u})\| \right) +\mu_0\|\textbf{F}\circ \textbf{G}_{m-1}\circ\cdots \circ \textbf{G}_1(\widetilde{u})-\textbf{G}_m \circ\cdots \circ \textbf{G}_1(\widetilde{u})\| \leq$$
$$ \sum_{p=1}^{m} \mu_{m-p}  \|\textbf{F}\circ \textbf{G}_{p-1}\circ\cdots \circ \textbf{G}_1(\widetilde{u})-\textbf{G}_p\circ\cdots \circ \textbf{G}_1(\widetilde{u})\|.$$ \hfill$\Box$\par

\noindent\textbf{Lemma 2.}  \textsl{Let $(X,\|\cdot\|)$ be a Banach space, $\textbf{F}: X \longrightarrow X$ and let $\mu_0=1$ and $\{\mu_n\}_{n \geq 1}$ be a sequence of nonnegative real numbers such that the series $\sum_{n \ge 1}\mu_n$ is convergent and  $\|\textbf{F}^n u_1-\textbf{F}^n u_2\|\leq \mu_n\|u_1-u_2\|$ for all $u_1, u_2 \in X$ and for all $n\in \mathbb{N}$.  Let  $u\in X$ be the unique fixed point of $\textbf{F}$, 
$\widetilde{u} \in X$, $\varepsilon >0$ and $m \in \mathbb{N}$ such that $\|u-\textbf{F}^m \widetilde{u}\| < \varepsilon/2$.  If $\textbf{G}_1,\cdots, \textbf{G}_m:X \longrightarrow X $ are operators such that 
$$ \|\textbf{F}^m(\widetilde{u})-\textbf{G}_m\circ...\circ \textbf{G}_1(\widetilde{u})\|<\varepsilon/2$$
then
$$\|u-\textbf{G}_m\circ...\circ \textbf{G}_1(\widetilde{u})\|<\varepsilon.$$}

 \bigskip
 
\textsl{Proof.}
By triangular inequality, we have that 
$$\|u-\textbf{G}_m\circ...\circ \textbf{G}_1(\widetilde{u})\| \leq \|u-\textbf{F}^m(\widetilde{u})\|+\|\textbf{F}^m(\widetilde{u})-\textbf{G}_m\circ...\circ \textbf{G}_1(\widetilde{u})\|    
< \frac{\varepsilon}{2}+ \frac{\varepsilon}{2} <\varepsilon.$$ \hfill$\Box$\par

\bigskip

Banach's Fixed-Point Theorem allows us to express the fixed point of $\textbf{F}$ as the limit of the sequence of functions $\{\textbf{F}^m (\widetilde{u})\}_{m\in \mathbb{N}},$ with $\widetilde{u} \in X$. Obviously, if it were possible to explicitly calculate, for each iteration, the expression $\textbf{F}^m (\widetilde{u})$, then for each $m$ we would have an approximation of the fixed point. But as a practical matter, such an explicit calculation is only possible in very particular cases. For this reason, in view of Lemma 2, for the problem at hand we will begin with an initial function $\widetilde{u} \in X$ and obtain successive $\textbf{G}_m\circ...\circ \textbf{G}_1(\widetilde{u})$ approximations of the fixed point $u$ of the operator $\textbf{F}$ following the scheme:

$$
\begin{array}{cccccccccc}
 \widetilde{u} &  &  &  &  &  &  &  &  &  \\
  \downarrow &  &  &  &  &  &  &  &  &  \\
  \textbf{F}(\widetilde{u}) & \approx & \textbf{G}_1(\widetilde{u}) &  &  &  &  &  &  &  \\
  \downarrow &  & \downarrow &  &  &  &  &  &  &  \\
  \textbf{F}^2(\widetilde{u})  & \approx & \textbf{G}_2\circ\textbf{G}_1(\widetilde{u}) &  &  &  &  &  \\
  \downarrow &  &  \downarrow &  &  &  &  &  \\
  \textbf{F}^3(\widetilde{u}) & \approx & \textbf{G}_3\circ\textbf{G}_2\circ \textbf{G}_1(\widetilde{u})&  &  &  \\
  \vdots &  &   \vdots &  &  \\
  \downarrow &  &   \downarrow &  &  \\
   \textbf{F}^m(\widetilde{u})  & \approx &
  \textbf{G}_m\circ...\circ \textbf{G}_1(\widetilde{u})\approx u \end{array}
$$

\bigskip

\noindent\textbf{Remark.} Note that:
\begin{enumerate}
\item[a)] By (\ref{1}) and the convergence of $\sum_{n \ge 1}\mu_n$, given $\varepsilon>0$, we can choose $m\in \mathbb{N}$ such that
$$\sum_{i=m}^{\infty}\mu_i<\frac{\varepsilon}{2 \|\textbf{F}\widetilde{u}-\widetilde{u}\|}$$
and that implies that $\|u-\textbf{F}^{m} \widetilde{u}\| < \varepsilon/2$.  .
\item[b)] If $ \|\textbf{F}\circ \textbf{G}_{p-1}\circ\cdots \circ \textbf{G}_1(\widetilde{u})-\textbf{G}_p\circ\cdots \circ \textbf{G}_1(\widetilde{u})\| <\dfrac{\varepsilon}{2 m \mu_{m-p}}, $ for $ p=1,\dots,m$
then by Lemma 1, $\|\textbf{F}^m(\widetilde{u})-\textbf{G}_m\circ...\circ \textbf{G}_1 (\widetilde{u})\|<\varepsilon/2.$ 
\end{enumerate}

Therefore theoretically, given $\varepsilon$ the number of $\textbf{G}_{i}$ needed to reach the accuracy can be determined by the above conditions.

 \section{Biorthogonal Systems}
 
 In the next section, we will apply the approximation provided by Lemma 2 to approximate the solution of NPVIDEs. In the construction of certain $\textbf{G}_p$ operators that verify the conditions of the aforementioned Lemma 2, we will use biorthogonal systems in Banach spaces. For completing the work, we collect the most important notations and properties of this tool.

Biorthogonal systems and in particular the Faber Schauder bases have been successfully combined with the fixed point theory (see \cite{brezis}, \cite{gelbaum}, \cite{jameson},  \cite{semadeni1}, and \cite{semadeni2})   to approximate the solution of differential, integral and integro-differential equations (see \cite{berenguer2009} -- \cite{berenguer1} and \cite{gamez}).

Let us start by recalling the  notion of  biorthogonal system of a Banach space. Let $X$ be a Banach space and $X^*$ its topological dual space. A system $\{\psi_n,f_n\}_{n\ge 1}$, where $\psi_n \in X$,
$f_n \in X^*$, and $f_n(\psi_m) =\delta_{n m}$ ($\delta$ is Kronecker's delta) is called a \textsl{biorthogonal system} in
$X$. We say that the system is a \textsl{fundamental} biorthogonal system if $\overline{\textrm{span}} \{\psi_n\}_{n\ge 1}=X$.

We will work with a particular type of fundamental biorthogonal system. Let us recall that a sequence $\{\psi_n\}_{n\ge 1}$  of elements of a Banach space $X$ is called a \textsl{Schauder basis} of $X$ if, for every $z\in X$, there is a unique sequence
$\{\lambda_n\}_{n\ge 1}$ of scalars such that $z=\sum_{n\ge 1} \lambda_n \psi_n$. A Schauder basis gives rise to the canonical sequence of (continuous and linear) finite dimensional \textsl{projections} $P_n:X\rightarrow X$, $P_n (\sum_{k\ge 1} \lambda_k \psi_k) = \sum_{k=1}^n \lambda_k \psi_k$; and the associated sequence of (conti\-nuous and linear) \textsl{coordinate functionals} $\{\psi_n^*\}_{n\ge 1}$ in  $X^*$ is defined  by
$\psi_n^*(\sum_{k\ge 1} \lambda_k \psi_k) = \lambda_n.$ Note that a Schauder basis is always a fundamental biorthogonal system, under the interpretation of the coordinate functionals as  biorthogonal functionals.

\bigskip
It follows from the Bayre category theorem (\cite{brezis}, Theorem 2.1) that for all $n \geq 1$, the biorthogonal functional $\psi_n^*$ and the projection $P_n$ are continuous; these notions
clearly imply that for all $z \in X$ there holds
\begin{equation}\label{2}
\begin{displaystyle}\lim_{n} \|P_n(z)-z\| \end{displaystyle}=0.
\end{equation}

 \section{Numerical solution of NPVIDEs}

The second aim of this work is to introduce a new numerical method  for solving a class of NPVIDEs. Specifically  let us consider the following nonlinear problem: given $a, g\in  C(\Omega)$, $K\in C(\Omega^2 \times \mathbb{R}$) and $u_0\in  C([\alpha,\alpha+\beta])$  with $\Omega=[\alpha,\alpha+\beta] \times
[0,T]$, find $u\in C(\Omega)$ with  $ \displaystyle \frac{\partial u}{\partial t}\in C(\Omega)$ such that:

\begin{equation}\label{3}
\left\{\begin{array}{l} \displaystyle
\frac{\partial u (x,t)}{\partial t}=a(x,t)u(x,t)+g(x,t)+\int_0^t\int_\alpha^x K(x,t,y,s,u(y,s)) dy\ ds    \\ \\
u(x,0)=u_0(x)

\end{array}
\right..
\end{equation}

We also assume that the kernel function $K$ satisfies a Lipschitz condition with respect to the last variable, with Lipschitz constant $M$.

Equation (\ref{3}) can be written in the abstract form: $\textbf{F}u=u$, where the operator $\textbf{F}:C(\Omega)\longrightarrow
C(\Omega)$ is defined for $(x,t)\in C(\Omega)$ and $u \in C(\Omega)$ as  

\begin{equation}\label{4}
(\textbf{F}u)(x,t):=u_0(x)+\int_0^t a(x,r)u(x,r)dr+\int_0^tg(x,r)dr +\int_0^t \int_0^r \int_{\alpha}^x K(x,r,y,s,u(y,s))dydsdr.
\end{equation}

It is clear that $u$ is a solution of (\ref{3}) if, and only if, it is a fixed point of the operator $\textbf{F}$.

For all $u_1,u_2 \in C(\Omega)$, $(x,t)\in \Omega$ and for all $n\geq 1$ it is easy to show by  mathematical induction that, if  $\begin{displaystyle} N= \end{displaystyle}\max_{(x,t) \in \Omega}|a(x,t)|$, then 
$$|(F^n u_1)(x,t)-(F^n u_2)(x,t)|\leq 
\left( \sum_{i=0}^n {n\choose i} \frac{N^{n-i}\ t^{n+i}\ (x-\alpha)^i}{(n+i)!\ i!} M^i  \right)  \|u_1-u_2\|. $$
As consequence,
$$\|\textbf{F}^n u_1-\textbf{F}^n u_2\| \leq \left( \sum_{i=0}^n {n\choose i} \frac{N^{n-i}\ T^{n+i}\ \beta^i}{(n+i)!\ i!} M^i  \right)  \|u_1-u_2\| =$$
$$ \frac{T^n}{n!}\left( \sum_{i=0}^n {n\choose i} \frac{N^{n-i}\ (T \beta M)^i}{(n+i)(n+i-1)\cdots (n+1) i!}   \right)  \|u_1-u_2\| \leq
\frac{T^n}{n!}\left( \sum_{i=0}^n {n\choose i} N^{n-i}\left(\frac{T \beta M}{n}\right)^i   \right)  \|u_1-u_2\|$$

Then for $n\geq 1$ \begin{equation}\label{Fn}
\|\textbf{F}^n u_1-\textbf{F}^n u_2\| \leq \mu_n \|u_1-u_2\| \end{equation}
with 
\begin{equation}\label{mun} \mu_n=\frac{T^n}{n!} \left(N+\frac{T \beta M}{n}\right)^n. \end{equation}

Since the series $\begin{displaystyle}\sum_{n\geq 1}\end{displaystyle} \mu_n$ is convergent, applying the Banach Fixed--Point Theorem gives us that $\textbf{F}$ has one and only one fixed point $u \in C(\Omega)$, and for each $\widetilde{u} \in C(\Omega)$ and $n\geq 1$,
\begin{equation}\label{5}
\|\textbf{F}^n \widetilde{u}-u\| \leq \sum_{i=n}^{\infty} \mu_i \|\textbf{F}\widetilde{u}-\widetilde{u}\|
\end{equation}
and thus $\begin{displaystyle}\lim_{ n\rightarrow\infty} \|\textbf{F}^n \widetilde{u}-u\| \end{displaystyle}=0.$

\subsection{Numerical method}

 Given $\varepsilon>0$, our goal in this section is to choose $ \textbf{G}_1, \textbf{G}_2,\cdots,\textbf{G}_m: C (\Omega) \longrightarrow C (\Omega) $, verifying the conditions of Lemma 2 to arrive at an  approximation of the solution to equation 
(\ref{3}).  We consider a Schauder basis $\{A_n\}_{n\geq 1}$ in $C(\Omega)$ and  the sequence of associated projections  $\{Q_n\}_{n\geq 1}$. In addition we consider $\{B_n\}_{n\geq 1}$ a Schauder basis in $C(\Omega^2)$, with its  associated projections  $\{R_n\}_{n\geq 1}$.  

The choice we propose is the following: 
\begin{equation}\label{6}
(\textbf{G}_p v)(x,t):=u_0(x)+\int_0^t Q_{n_p^2} (a \cdot v)(x,r)dr+\int_0^t g(x,r)dr +\int_0^t \int_0^r \int_{\alpha}^x R_{n_p^4}(  L_0(v))(x,r,y,s)dydsdr
\end{equation} 

\noindent where $n_p\in \mathbb{N}$,  $v \in C(\Omega)$  and $L_0:C(\Omega) \rightarrow C(\Omega^2)$,
\begin{equation} \label{7}
L_{0}(v)(x,t,y,s):=K(x,t,y,s,v(y,s)), \ \ \ \ \ \ (x,t,y,s) \in \Omega^2.
\end{equation}

Our objective is to justify that we can choose $n_1,\cdots,n_m \in \mathbb{N} $, so that the operators $ \textbf{G}_1,\cdots,\textbf{G}_m $ can be used to obtain an approximation of the unique solution of the equation (\ref{3}):
$$ u
\approx \textbf{G}_m\circ ...\circ \textbf{G}_1(\widetilde{u}) \ \ \ \ \widetilde{u}\in C(\Omega)$$

To this end, for $p=\{1,\dots,m\}$ we will denote $L_p:C(\Omega) \rightarrow C(\Omega^2)$,
\begin{equation} \label{10}L_{p}(v)(x,t,y,s):=K(x,t,y,s,\textbf{G}_p\circ \textbf{G}_{p-1}\circ\cdots \circ\textbf{G}_1(v)(y,s)), \ \ \ \ \ \ (x,t,y,s) \in \Omega^2.\end{equation} 

\noindent\textbf{Lemma 3.}

 \textsl{With the previous notation, given $\widetilde{u}\in C(\Omega)$ it is satisfied that   }
$$\|\textbf{F}^m (\widetilde{u})-\textbf{G}_m\circ ...\circ \textbf{G}_1(\widetilde{u})\| \leq $$ 
\begin{equation}\label{8}\sum_{p=1}^m \mu_{m-p} \left(  \| a\cdot( \textbf{G}_{p-1}\circ\cdots \circ \textbf{G}_1(\widetilde{u}))- Q_{n_p^2}(a\cdot (\textbf{G}_{p-1}\circ ...\circ \textbf{G}_1(\widetilde{u})))\|T+\|L_{p-1}(\widetilde{u})- R_{n_p^4}(L_{p-1}(\widetilde{u}))\|\beta T^2\right)
\end{equation}
\noindent with 
$$
\mu_{m-p}=\frac{T^{m-p}}{(m-p)!} \left(N+\frac{T \beta M}{m-p}\right)^{m-p} \textrm{ for } p\in\{1,\dots,m-1\} \textrm{ and } \mu_0=1.
$$

\textsl{Proof}

From (\ref{Fn}), (\ref{mun}) and Lemma 1, $$ \|\textbf{F}^m(\widetilde{u})-\textbf{G}_m\circ...\circ \textbf{G}_1(\widetilde{u})\|\leq \sum_{p=1}^{m} \mu_{m-p}  \|\textbf{F}\circ \textbf{G}_{p-1}\circ\cdots \circ \textbf{G}_1(\widetilde{u})-\textbf{G}_p\circ\textbf{G}_{p-1}\circ\cdots \circ \textbf{G}_1(\widetilde{u})\|, 
 $$
 where $ \mu_{m-p}$ is given by (\ref{mun}) for $p\in\{1,\dots,m-1\}$ and $\mu_0=1$.
 
On the other hand for $p\in\{1,\dots,m\}$,
$$ |\textbf{F}\circ \textbf{G}_{p-1}\circ\cdots \circ \textbf{G}_1(\widetilde{u})(x,t)-\textbf{G}_p\circ\textbf{G}_{p-1}\circ \cdots \circ \textbf{G}_1(\widetilde{u})(x,t)| $$
$$\leq \int_0^t\left| a\cdot (\textbf{G}_{p-1}\circ\cdots \circ \textbf{G}_1(\widetilde{u}))(x,r)-Q_{n_p^2}(a\cdot(\textbf{G}_{p-1}\circ ...\circ \textbf{G}_1(\widetilde{u})))(x,r)dr\right|$$
$$+\int_0^t \int_0^r \int_{\alpha}^x \left|L_{p-1}(\widetilde{u})(x,r,y,s)- R_{n_p^4}(L_{p-1}(\widetilde{u}))(x,r,y,s)  \right|dydsdr$$
$$\leq \| a\cdot( \textbf{G}_{p-1}\circ\cdots \circ \textbf{G}_1(\widetilde{u}))- Q_{n_p^2}(a\cdot (\textbf{G}_{p-1}\circ ...\circ \textbf{G}_1(\widetilde{u})))\|T+\|L_{p-1}(\widetilde{u})- R_{n_p^4}(L_{p-1}(\widetilde{u}))\|\beta T^2.$$
Therefore
$$ \|\textbf{F}^m(\widetilde{u})-\textbf{G}_m\circ...\circ \textbf{G}_1(\widetilde{u})\|\leq 
 $$
  $$\leq \sum_{p=1}^m \mu_{m-p}  \left( \| a\cdot( \textbf{G}_{p-1}\circ\cdots \circ \textbf{G}_1(\widetilde{u}))- Q_{n_p^2}(a\cdot (\textbf{G}_{p-1}\circ ...\circ \textbf{G}_1(\widetilde{u})))\|T+\|L_{p-1}(\widetilde{u})- R_{n_p^4}(L_{p-1}(\widetilde{u}))\|\beta T^2.\right)$$
\hfill$\Box$\par

\bigskip

\noindent\textbf{Theorem 4.} 

\textsl{With the previous notation, let $u$ be the unique solution of (\ref{3}). Then, for each $\varepsilon>0$ and each $\widetilde{u} \in C(\Omega)$, there exist $m \geq 1$ and $n_1,\cdots,n_m \geq 1$ such that $\|u-\textbf{G}_m\circ ...\circ \textbf{G}_1(\widetilde{u})\|<\varepsilon.$}

\textsl{Proof}

By (\ref{limite}), for $\varepsilon >0$, there exists $m\in \mathbb{N}$ such that $\|u-\textbf{F}^m \widetilde{u}\| < \varepsilon/2$. For that $m$,  in view of the convergence property (\ref{2}) and Lema 3, we can find
$n_1,\cdots,n_m \geq 1$ and therefore $\textbf{G}_1,\cdots,\textbf{G}_m$, such that
$$\|\textbf{F}^m(\widetilde{u})-\textbf{G}_m\circ...\circ \textbf{G}_1(\widetilde{u})\|\leq $$
$$\sum_{p=1}^m \mu_{m-p}  \left( \| a\cdot( \textbf{G}_{p-1}\circ\cdots \circ \textbf{G}_1(\widetilde{u}))- Q_{n_p^2}(a\cdot (\textbf{G}_{p-1}\circ ...\circ \textbf{G}_1(\widetilde{u})))\|T+\|L_{p-1}(\widetilde{u})- R_{n_p^4}(L_{p-1}(\widetilde{u}))\|\beta T^2\right)<\varepsilon/2.$$

Then for Lemma 2,
$$\|u-\textbf{G}_m\circ ...\circ \textbf{G}_1(\widetilde{u})\|<\varepsilon.$$

\hfill$\Box$\par

 \bigskip
 \subsection{Numerical examples}
 
 This section contais some examples and their numerical results. In order to obtain concrete numerical approximations of the solution of NPVIDE, we shall fix concrete bivariate and multivariate Schauder bases in $C(\Omega)$ and $C(\Omega^2)$, respectively, called usual Schauder bases, derived from the usual Schauder bases in $C([\alpha,\alpha+\beta])$ and $C([0,T])$. 

Let $\{x_n\}_{n\geq 1}$ be a sequence of distinct points from $[\alpha,\alpha+\beta]$ such that $x_1=\alpha$ and $x_2=\alpha+\beta$, the usual Schauder basis $\{\varphi_n\}_{n\geq 1}$ in $C([\alpha,\alpha+\beta])$ being defined by $\varphi_1(x):=1 \ x \in [\alpha,\alpha+\beta]$; and for $n\geq 2$, $\varphi_n$ is the piecewise linear continuous function on $[\alpha,\alpha+\beta]$ with nodes at $\{x_i: 1 \leq i \leq n\}$ uniquely determined by the relation $\varphi_n(x_n)=1$ and for $k<n, \ \varphi_n(x_k)=0$. Analogously, for a dense sequence $\{t_n\}_{n\geq 1}$ of distinct points in $[0,T]$, with $t_1=0$ and $t_2=T$, we will denote as $\{\widehat{\varphi}_n\}_{n \geq 1}$
the usual Schauder basis in $C([0,T])$ (see \cite{gelbaum} and \cite{semadeni2}).

We consider the bijective mappings $\tau:\mathbb{N} \longrightarrow
\mathbb{N}\times \mathbb{N}$, and 
$\phi:\mathbb{N} \longrightarrow
\mathbb{N}\times \mathbb{N}\times \mathbb{N}\times \mathbb{N}$, defined respectively by $\tau(n):=(\tau_1(n),\tau_2(n))$ and $\phi(n):=
(\tau_1(\tau_1(n)), \tau_2(\tau_1(n)),\tau_1(\tau_2(n)),\tau_2(\tau_2(n))),$
where

\noindent $$\tau(n):=\left\{
\begin{array}{cl}
  (\sqrt{n},\sqrt{n}) &  \hbox{if} \quad [\sqrt{n}]=\sqrt{n} \\

(n-[\sqrt{n}]^2,[\sqrt{n}]+1) & \hbox{if} \quad
 0<n-[\sqrt{n}]^2\le [\sqrt{n}] \\

 ([\sqrt{n}]+1,n-[\sqrt{n}]^2-[\sqrt{n}]) &  \hbox{if} \quad
  [\sqrt{n}]<
  n-[\sqrt{n}]^2
\end{array}
\right.,$$

\noindent where $[q]$ denotes the integer part of $q$.

From the Schauder bases $\{\varphi_n\}_{n \geq 1}$ and 
$\{\widehat{\varphi}_n\}_{n \geq 1}$, we can build another usual Schauder basis $\{A_n\}_{n \geq 1}$ of $C(\Omega)$. It is sufficient to consider 
$$A_n(x,t):=\varphi_i (x) \widehat{\varphi}_j (t), \ \  (x,t) \in \Omega$$
with $\tau(n)=(i,j)$ (see \cite{berenguer4}). The usual Schauder basis $\{B_n\}_{n \geq 1}$ of $C(\Omega^2)$ (see \cite{berenguer2}) is defined as:
$$B_n(x,t,y,s):=\varphi_i (x) \widehat{\varphi}_j (t) \varphi_k (y)  \widehat{\varphi}_l (s), \ \ \
x,y \in [\alpha,\alpha+\beta], \ t,s \in [0,T],$$
whenever $\phi(n)=(i,j,k,l)$.

For the usual Schauder bases chosen, we could estimate the convergence rate of the sequence of projections taking into account the properties of such bases and the Mean-Value Theorem  in a similar way to \cite{berenguer1} and \cite{berenguer2}.

In the following examples we consider $m=2$ and $m=3$. The absolute errors obtained by present method ($e_2(x,t)=\vert G_2\circ G_1(\widetilde{u})(x,t)-u(x,t)\vert $ and  $e_3(x,t)=\vert G_3\circ G_2\circ G_1(\widetilde{u})(x,t)-u(x,t)\vert $) are shown in tables with $\widetilde{u}(x,t)=u_0(x)$ and $n_p=3$ for $1\leq p\leq m$. All numerical results are obtained by running a program written in Mathematica software.  
 
 \noindent\textbf{Example 1.} 
 
In this example, we consider $a(x,t)= t\sin(x)$, $K(x,t,y,s,u)= y s u$, $g(x,t)=x-\dfrac{t^3 x^3}{9}-t^2 x \sin(x)$, $u_0(x)=0$ and where  $u(x,t)=x t$ is the exact solution. The absolute errors at specified points are reported in Table 1.  Also the last row includes the errors  with the collocation method with bivariate classical polinomial interpolation.

\begin{center}

\small{\textsc{Table 1. Numerical results for Example 1.}}

\noindent \begin{tabular}{|c|c|c|c|c|c|c|c|c|c|c|c| }\hline 
  $(x,t)$& $(0,0)$ & $(0.1,0.1)$& $(0.2,0.2)$ & $(0.4,0.4)$  &$(0.6,0.6)$&  $(0.8,0.8)$ &  $(0.9,0.9)$&  $(1,1)$  \\ \hline
 $e_2(x,t)$& $0$ & $4.76\times 10^{-4}$ & $3.73\times 10^{-3}$ & $2.75\times 10^{-2}$& $7.79\times 10^{-2}$ & 
  $1.28\times 10^{-1}$ & $1.38\times 10^{-1}$ & $2.01\times 10^{-1}$
   
      \\ \hline  
  
 $e_3(x,t)$& $0$ & $3.33\times 10^{-6}$ & $1.06\times 10^{-4}$ & $3.36\times 10^{-3}$& $2.51\times 10^{-2}$ & 
  $1.03\times 10^{-1}$ & $1.08\times 10^{-1}$ & $1.25\times 10^{-1}$\\ \hline

 $\textrm{Collocation}$& $0$ & $4.99\times 10^{-2}$ & $9.85\times 10^{-2}$ & $1.96\times 10^{-1}$& $3.03\times 10^{-1}$ & 
  $4.29\times 10^{-1}$ & $5.01\times 10^{-1}$ & $5.78\times 10^{-1}$\\ \hline

\end{tabular}
\end{center}

 \noindent\textbf{Example 2.}  Now we consider $a(x,t)= 0$, $K(x,t,y,s,u)= u t^2 $, $g(x,t)=x \cos(t)-t^2 x^2 \sin^2\left(\dfrac{t}{2}\right)$, $u_0(x)=0$ and where  $u(x,t)=x \sin(t)$ is the exact solution.  Its numerical
results are given in Table 2.

\begin{center}

\small{\textsc{Table 2. Numerical results for Example 2.}}

\noindent \begin{tabular}{|c|c|c|c|c|c|c|c|c|c|c|c| }\hline 
  $(x,t)$& $(0,0)$ & $(0.1,0.1)$& $(0.2,0.2)$ & $(0.4,0.4)$  &$(0.6,0.6)$&  $(0.8,0.8)$ &  $(0.9,0.9)$&  $(1,1)$  \\ \hline
 $e_2(x,t)$& $0$ & $4.78\times 10^{-4}$ & $3.83\times 10^{-3}$ & $1.02\times 10^{-2}$& $3.01\times 10^{-2}$ & 
  $2.34\times 10^{-1}$ & $3.27\times 10^{-1}$ & $4.41\times 10^{-1}$ \\ \hline

  $e_3(x,t)$& $0$ & $4.99\times 10^{-9}$ & $6.38\times 10^{-7}$ & $8.11\times 10^{-5}$& $1.36\times 10^{-3}$ &  $1.00\times 10^{-2}$ &  $2.27\times 10^{-2}$ &   $4.7\times 10^{-2}$  \\ \hline

\end{tabular}
\end{center}

\noindent\textbf{Example 3.}   We consider $a(x,t)= x \sin(t)$, $K(x,t,y,s,u)= x t^2 u $, $g(x,t)=-\dfrac{t^4 x }{2}+\dfrac{1}{2} t^4 x \cos(x)+\sin(x)-t x \sin(t) \sin(x)$, $u_0(x)=0$ and where  $u(x,t)=t \sin(x)$ is the exact solution.  Its numerical results are shown in Table 3. Also last row includes the errors  with the collocation method.

\begin{center}

\small{\textsc{Table 3. Numerical results for Example 3.}}

\noindent \begin{tabular}{|c|c|c|c|c|c|c|c|c|c|c|c| }\hline 
  $(x,t)$& $(0,0)$ & $(0.1,0.1)$& $(0.2,0.2)$ & $(0.4,0.4)$  &$(0.6,0.6)$&  $(0.8,0.8)$ &  $(0.9,0.9)$&  $(1,1)$  \\ \hline
  
$e_2(x,t)$& $0$ & $3.56\times 10^{-4}$ & $3.57\times 10^{-3}$ & $2.61\times 10^{-2}$& $6.95\times 10^{-2}$ & 
  $1.07\times 10^{-1}$ & $1.77\times 10^{-1}$ & $2.99\times 10^{-1}$ \\ \hline

  $e_3(x,t)$& $0$ & $3.32\times 10^{-6}$ & $1.05\times 10^{-4}$ & $3.30\times 10^{-3}$& $2.45\times 10^{-2}$ &  $9.97\times 10^{-2}$ & $9.8\times 10^{-2}$ &$7.37\times 10^{-2}$ \\ \hline

 $\textrm{Collocation}$&  $0$ & $6.44\times 10^{-2}$ & $2.20\times 10^{-1}$ & $1.60\times 10^{-1}$& $2.52\times 10^{-1}$ & 
  $5.34\times 10^{-1}$ & $4.74\times 10^{-1}$ & $3.35\times 10^{-1}$\\ \hline

\end{tabular}
\end{center}

 \noindent\textbf{Example 4.}   In the last example, we consider $a(x,t)= x^2$, $K(x,t,y,s,u)= t \sin(u) $, $g(x,t)=2-x^2(2t+x)-t \sin(t)(\cos(t)-\cos(t+x))$, $u_0(x)=x$ and where  $u(x,t)=x+2t$ is the exact solution.  The numerical results  are given in Table 4.

\begin{center}

\small{\textsc{Table 4. Numerical results for Example 4.}}

\noindent \begin{tabular}{|c|c|c|c|c|c|c|c|c|c|c|c| }\hline 
  $(x,t)$& $(0,0)$ & $(0.1,0.1)$& $(0.2,0.2)$ & $(0.4,0.4)$  &$(0.6,0.6)$&  $(0.8,0.8)$ &  $(0.9,0.9)$&  $(1,1)$  \\ \hline
  
$e_2(x,t)$& $0$ & $1.81\times 10^{-2}$ & $2.12\times 10^{-2}$ & $3.03\times 10^{-2}$& $5.20\times 10^{-2}$ & 
  $1.57\times 10^{-1}$ & $2.29\times 10^{-1}$ & $1.93\times 10^{-1}$ \\ \hline
  
  $e_3(x,t)$& $0$ & $1.47\times 10^{-3}$ & $5.67\times 10^{-3}$&$1.87\times 10^{-2}$& $2.73\times 10^{-2}$ &  $5.33\times 10^{-2}$ & $7.48\times 10^{-2}$ & $1.27\times 10^{-1}$ \\ \hline

\end{tabular}
\end{center}
 
\section*{Acknowledgements}

\noindent  The authors would like to thank the referee for her/his careful reading.

\noindent Research partially supported by project MTM2016-80676-P (AEI/Feder, UE), by Junta de Andaluc\'{i}a
Grant FQM359 and by E.T.S. Ingenier\'{i}a de Edificaci\'{o}n of
the University of Granada (Spain).

%
%

\end{document}